\documentclass[12pt]{amsart}

\usepackage[lite]{amsrefs}

\usepackage{verbatim}
\usepackage{amssymb}
\usepackage{upref}

\newtheorem{thm}{Theorem}[section]
\newtheorem{cor}[thm]{Corollary}
\newtheorem{lem}[thm]{Lemma}
\newtheorem{prop}[thm]{Proposition}

\theoremstyle{definition}
\newtheorem{defn}[thm]{Definition}

\newtheorem{stand}[thm]{Standing Hypothesis}

\newtheorem{rem}[thm]{Remark}

\numberwithin{equation}{section}

\newcommand{\secref}[1]{Section~\textup{\ref{#1}}}
\newcommand{\thmref}[1]{Theorem~\textup{\ref{#1}}}
\newcommand{\corref}[1]{Corollary~\textup{\ref{#1}}}
\newcommand{\lemref}[1]{Lemma~\textup{\ref{#1}}}
\newcommand{\propref}[1]{Proposition~\textup{\ref{#1}}}
\newcommand{\defnref}[1]{Definition~\textup{\ref{#1}}}
\newcommand{\remref}[1]{Remark~\textup{\ref{#1}}}

\newcommand{\righttext}[1]{\qquad\text{#1 }}

 \DeclareMathOperator*{\spn}{span}

\newcommand{\cc}[1]{\mathcal{#1}}
\newcommand{\bb}[1]{\mathbb{#1}}

\newcommand{\R}{\bb R}
\newcommand{\C}{\bb C}

\newcommand{\A}{\cc A}
\newcommand{\B}{\cc B}
\renewcommand{\L}{\cc L}
\renewcommand{\P}{\cc P}

\newcommand{\Chi}{\raisebox{2pt}{\ensuremath{\chi}}}

\newcommand{\SL}{\mathrm{SL}}
\newcommand{\NN}{\cc N}
\newcommand{\MM}{\cc M}
\newcommand{\RR}{\cc R}

\newcommand{\<}{\langle}
\renewcommand{\>}{\rangle}
\renewcommand{\iff}{\Leftrightarrow}
\newcommand{\bsl}{\setminus}

\newcommand{\inv}{^{-1}}

\newcommand{\norm}[1]{\lVert{#1}\rVert}
\newcommand{\Norm}[1]{\left\lVert{#1}\right\rVert}

\newcommand{\ie}{\emph{i.e.}}
\newcommand{\eg}{\emph{e.g.}}
\newcommand{\cf}{\emph{cf.}}
\newcommand{\st}{\emph{s.t. }}

\begin{document}

\title[Compact open subgroups and multiplier Hopf $^*$-algebras]
{Groups with compact open subgroups and multiplier Hopf $^*$-algebras}

\author[Landstad]{Magnus B. Landstad}
\address{Department of Mathematical Sciences\\
Norwegian University of Science and Technology\\
NO-7491 Trondheim, Norway} \email{magnusla@math.ntnu.no}

\author[Van Daele]{A. Van Daele}
\address{Department of Mathematics\\
K.U.~Leuven, Celestijnenlaan 200B, BE-3001 Heverlee, Belgium}
\email{Alfons.VanDaele@wis.kuleuven.be}

\subjclass[2000]{ 22D05;  46L05}

\keywords{Totally disconnected
groups,  group $C^*$-algebras, multiplier Hopf algebras.}

\begin{abstract}
For a locally compact group $G$  we look at the group algebras
$C_0(G)$ and $C_r^*(G)$,
and we let $f\in C_0(G)$ act on $L^2(G)$ by the multiplication operator $M(f)$.
We show among other things that the following properties are equivalent:

1. $G$ has a compact open subgroup.

2.  One of the $C^*$-algebras has a dense multiplier Hopf $^*$-sub\-alg\-e\-bra (which turns out to be unique).

3.  There are non-zero elements $a\in C_r^*(G)$ and  $f\in C_0(G)$ such that $aM(f)$ has finite rank.

4.  There are non-zero elements $a\in C_r^*(G)$ and  $f\in C_0(G)$ such that
$aM(f)=M(f)a$.

If $G$ is abelian, these properties are equivalent to:

5.  There is a non-zero continuous function
 with the property that both $f$ and $\widehat f$ have compact support.
\end{abstract}

\maketitle


\section*{Introduction}
\label{sec:Intro}

The background for this article is the observation in \cite{KLQ}*{Section~3}
 that the
algebra $C_c^\infty(G)$ of {\it smooth} functions on a totally disconnected locally
compact group $G$
is a {\it multiplier Hopf $^*$-algebra} as defined in \cite{VD-MH}. Therefore it
is natural to classify all
multiplier Hopf $^*$-algebras which are commutative or co-commutative; this is the
same as answering the following question: When does
$C_0(G)$ or $C_r^*(G)$
have a dense
multiplier Hopf $^*$-algebra? It is well known that the answer is yes if $G$ is
compact or discrete and the main results of \secref{sec:c-null-G} and
\secref{sec:c*G}  are
 that  in general the answer is yes if and only if $G$ has a
compact open subgroup.
If so, such a multiplier Hopf $^*$-algebra is unique and can be described explicitly as the algebra of polynomial functions on $G$.

We believe that multiplier Hopf $^*$-algebras in general are the right framework for
studying
totally disconnected locally compact quantum groups (what ever that is) and it is
therefore natural to first give a complete account for $C_0(G)$ and $C_r^*(G)$.

$C_0(G)$  is treated in \secref{sec:c-null-G}, the main tool is in
\corref{cor:mult.C-null-G} where it is shown that the existence of functions satisfying
some
algebraic relations is equivalent to
the existence of a compact open subgroup.

To get the same results for $C_r^*(G)$ is a little more
tricky.
For  $C_0(G)$ one can quickly show that elements of a multiplier Hopf
$^*$-subalgebra must have compact support and are therefore integrable.
It is not so easy to show that elements of a multiplier Hopf $^*$-subalgebra of
$C_r^*(G)$ automatically are integrable with respect to the Haar-Plancherel weight.
However, when this is proved, the results for $C_r^*(G)$ follow from those of
$C_0(G)$.
In fact, we show that if $\A$ is the unique dense
multiplier Hopf $^*$-algebra of $C_0(G)$ and $L$ denotes the left regular  representation, then
$\{L(f)\mid f\in\A\}$ is the unique dense
multiplier Hopf $^*$-algebra of $C_r^*(G)$.

For $C^*(G)$ the situation is different.
Here the existence of a multiplier Hopf $^*$-subalgebra does not imply
that $G$ has a
compact open subgroup. However,  the corresponding uniqueness result is true.

The algebras $C_0(G)$ and $C_r^*(G)$ are dual as locally compact quantum groups
 and we shall also see that
many properties of this duality are equivalent with the existence of a compact open subgroup.

It is well known that if $a\in C_r^*(G)$, $f\in C_0(G)$ and $M(f)$ is the the
corresponding multiplication operator on $L^2(G)$, then  $aM(f)$ is compact.
Since the finite rank operators are
dense in the algebra of compact operators it is natural to ask when $aM(f)\neq 0$ is
of finite rank.
We show this is possible if and only if $G$ has a compact open subgroup.

It is a consequence of the Heisenberg relations that if $G=\R^n$ then $a$ and $M(f)$
as above never commute unless one of them is zero. We show that in general $aM(f)=M(f)a\neq0$
is possible if and only if $G$ has a compact open subgroup.

Finally, as a bonus for the patient reader we look at the case where $G$ is abelian and ask whether one can have $f\in
C_c(G)$,
$f\neq0$ and $\widehat f\in C_c(\widehat G)$. The answer should not surprise.

Some of the results here are probably folklore and known to those who
have worked with representations of $p$-adic groups.

{\bf Acknowledgements.} We would like to thank the referee for valuable comments which in particular has shortened and improved the presentation in Section 1.
This research has been going on over a long
period, starting at The Centre for Advanced Study (CAS) in Oslo, and
continuing at the K.U.~Leuven and NTNU. We thank  for their hospitality and financial support, in addition we also have received support from The Research
Council of Norway (NFR) and The Research Council of  the
K.U.~Leuven.

 \section{Preliminaries}\label{sec:prelim} We start with fixing our notation regarding function spaces on $G$.

\begin{defn}
\label{defn:reg.repr}
As usual $C_0(G)$ are the continuous complex functions on $G$
vanishing at $\infty$ and $L^p(G)$  is defined with respect to a fixed left
Haar measure $\mu$.
We will often write just
\[
\int fd\mu=\int f(x)\,dx.
\]
The left and right action of $G$ on such functions is given by
\begin{equation}
_xf(y)=f(x\inv y)\qquad f_x(y)=f(yx).
\end{equation}
However, for functions in $L^2(G)$ we will instead use the left and right regular
representations given by
\begin{equation}
\label{leftrightreg}
L_xf(y)=f(x\inv y)\qquad   R_xf(y)=\Delta_G(x)^{1/2}f(yx)
\end{equation}
where $\Delta_G$ is the modular function on $G$.
\end{defn}

Many of the arguments are based on the
following Lemma.

\begin{lem}
\label{lem:mult.Uopen} Suppose $G$ is a locally compact group with a
continuous action $\alpha$ by isometries on a Banach space $A$.
For a fixed non-zero vector $a\in A$ and any subset $U\subset G$,
denote by $F_U(a)$ the linear span of the set
$\{ \alpha_x(a) \mid x\in U  \}$.

\begin{enumerate}
\item
Suppose that   $F_U(a)$ is finite dimensional
for a neighborhood $U$ of $e$.
 Then $G$ has an open subgroup $H$ \st $F_H(a)$ is also  finite dimensional.
\item
The same conclusion holds if we have a non-negligible set $C$ \st $F_C(a)$ is finite dimensional.
\item If one further assumes that the functions
$x\mapsto \<\alpha_x(a), \phi\>$
are in $C_0(G)$ for $a\in A$, $\phi\in A^*$, then the subgroup $H$ is also compact.
\end{enumerate}
\end{lem}

\begin{proof}
Clearly $V\subset U$ implies that $F_V(a)\subset F_U(a)$. To prove (i), take a neighborhood $U$ \st $F_U(a)$ has minimal, positive dimension. Take a neighborhood $V$ of $e$ \st $V=V\inv$
and $V^2\subset U$. Then $F_V(a)=F_U(a)$ is invariant by the open subgroup $H$ generated by $V$ and therefore $F_V(a)=F_H(a)$.

For (ii), take a measurable set $C$ with finite Haar measure $\mu(C)>0$ \st $F_C(a)$ has minimal, positive dimension.
Clearly
$U=\{y\in G\mid \mu(y\inv C\cap C)>0$ is a neighborhood of $e$ which satisfies
    \[
    y\in U \implies F_C(a)=F_{y\inv C\cap C}(a).
\]

 Therefore, the non-trivial linear space $F_C(a)$ is invariant by the open subgroup $H$  generated by $U$ and also here  $F_C(a)=F_H(a)$ is finite dimensional.

As for (iii),  let
$\{b_i\}$ be a basis for $F_H(a)$,  let $\phi_j\in A^*$ \st
$\<b_i, \phi_j\>=\delta_{ij}$ and define
$\psi_{ji}(z)=\<\alpha_z(b_i), \phi_j\>$.
Then $\psi_{ji}\in C_0(G)$ and
for $y, z \in H$:
\begin{align*}
\alpha_z(b_i)&=\sum_k \psi_{ki}(z)b_k\\
\psi_{ji}(yz)&=\sum_k \psi_{ki}(z)\<\alpha_y(b_k), \phi_j\>\\
             &=\sum_k \psi_{jk}(y)\psi_{ki}(z).
\end{align*}
So $1=\psi_{11}(yy\inv)=\sum_k \psi_{1k}(y)\psi_{k1}(y\inv) $ is
constant on $H$ and in $C_0(G)$, so $H$ must be compact.
\end{proof}
\begin{cor}
\label{cor:mult.C-null-G} If  $f,g, f_i,
g_i\in C_0(G)$ are non-zero functions \st
\begin{equation}
\label{mult.C-null-G} f(xy)g(y)=\sum_1^n f_i(x)g_i(y) \righttext{for all}
x,y \in G,
\end{equation}
then $G$ has a compact open subgroup $H$ and there are functions
$f'_j\in C_0(G)$ and $g'_j\in C(H)$ \st
\begin{equation}
\label{mult.C-null-Gb} f(xy)=\sum_1^n f'_j(x)g_j'(y) \righttext{for all} x
\in G, y\in H.
\end{equation}
\end{cor}
\begin{proof}
Pick a non-empty open subset $U$ of $G$ with $g(y)\neq 0$ for
$y\in U$, and we may assume  $e\in U$. We then have
\begin{equation*}
f_y=\sum  h_i(y)f_i \quad\text{where}\quad
h_i(y)={g_i}(y)/{g}(y) \quad\text{for} \quad y\in  U ,
\end{equation*}
so the result follows from \lemref{lem:mult.Uopen} by taking $A=C_0(G)$
and the  action given by $\alpha_y(f)= f_y$.
\end{proof}
\begin{lem}
\label{lem:comp.supp} Suppose that we have a continuous action of $G$
as in \lemref{lem:mult.Uopen}, that we have a finite set of non-zero
elements $a, a_i\in A$, and functions $g, g_i$ \st
\begin{equation}
\label{comp.supp} g(y)\alpha_y(a)=\sum g_i(y)a_i \righttext{for}y\in
G.
\end{equation}
Then $g$ has compact support.
\end{lem}
\begin{proof}
We may assume that $\{a_i\}$ are linearly independent, so
pick  $\nu_i\in A^*$ \st
$\<a_i,\nu_j\>=\delta_{ij}$.
Then
\begin{align*} g_i(y)&=g(y) \<\alpha_y(a),\nu_i\>\\
g(y)\alpha_y(a)&=g(y) \sum \<\alpha_y(a),\nu_i\>a_i\\
g(y)a&=g(y) \sum \<\alpha_y(a),\nu_i\>\alpha_{y\inv}(a_i).
\end{align*}
Pick $\nu_0\in A^*$ with $\<a,\nu_0\>=1$, then
\[
g(y)=g(y) \sum
\<\alpha_y(a),\nu_i\>\<\alpha_{y\inv}(a_i),\nu_0\>=g(y)c(y)
\]
 with
$c\in C_0(G)$,  which is possible only if $g$ has compact support.
\end{proof}

\begin{cor}
\label{cor:comp.supp} With the  assumptions  in
\corref{cor:mult.C-null-G}, $g\in C_c(G)$.
\end{cor}

We shall also need the following:
\begin{lem}
\label{lem:finite.support} Suppose $G$ contains a subgroup $H$, that
$A$ is a vector space and that we have functions $g,g_i:G\mapsto
A$ and $f_i:G\mapsto \C$ \st
\begin{equation}
\label{finite.support} \Chi_H(xy)g(y)=\sum_1^n f_i(x)g_i(y)
\righttext{for all} x,y \in G.
\end{equation}
Then $g$ has finite support in $G/H$, \ie\ there is a finite set
$F$ \st $g(y)=0$ for $y\notin FH$ .
\end{lem}
\begin{proof}
We may assume that the set $\{f_i\}$ is linearily independent, so
by \cite{HR2}*{(28.14)} there is a finite set $F_1=\{x_i\}$ \st the matrix  $\{f_i(x_j)\}$
is invertible. Then $y\notin F_1\inv H\Rightarrow g_i(y)=0$ and
therefore also $g(y)=0$. So we can take $F=F_1\inv$.
\end{proof}
\begin{cor}
\label{cor:multL2} Suppose the functions $f,g, f_i,
g_i\in L^2(G)$ satisfy
\begin{equation}
\label{multL2} f(x)g(x\inv y)=\sum_1^n f_i(x)g_i(y) \righttext{for almost all}
x,y \in G.
\end{equation}
Then $G$ has a compact open subgroup $H$ and there are functions
$f'_j\in C(H)$ and $g'_j\in L^2(G)$ \st
\begin{equation}
\label{multL2b} g(x\inv y)=\sum_1^m f'_j(x)g_j'(y) \righttext{ for }
x\in H, y\in G.
\end{equation}
\end{cor}
\begin{proof}
Pick a set $C$ such that $0<\mu(C)<\infty$ and $f(x)\neq 0$ for $x\in C$.
Then divide by $f(x)$ and apply
part (ii) of
\lemref{lem:mult.Uopen} with $A=L^2(G)$
and the action given by $L_x$.
\end{proof}
\begin{rem}
\label{cor:W-operator} The Kac-Takesaki operator on $ L^2(G\times G)$ is defined by
$
Wf(x,y)=f(x,x\inv y)$. Therefore
(\ref{multL2b}) is equivalent to
$W(f\otimes g)=\sum_1^n f_i\otimes g_i$.
The reader may want to compare this with \cite{BS}*{Definition~1.8}.
\end{rem}
\begin{rem}
Our main results so far have been that the existence of certain
functions satisfying \eqref{mult.C-null-G} is possible only if $G$ has a
compact open subgroup. Note, however, that this conclusion is
possible only with some restriction on the functions involved.
On $\R$ or matrix groups like $GL(n,\C)$ one clearly has unbounded functions
satisfying \eqref{mult.C-null-G}, but there are no compact open subgroups.
\end{rem}
\begin{rem}
As we shall see later, the conditions studied here are in fact
equivalent to the existence of a compact open subgroup $H$. For
the opposite implication, just take $f=g=\Chi_H$.
\end{rem}

\section{Multiplier Hopf $^*$-Algebras}\label{sec:mult.Hopf}
Multiplier Hopf $^*$-algebras were introduced in \cite{VDZ}, in this section we shall
recall some of the main definitions and results  and refer to
\cite{VDZ} or \cite{LvD} for
more precise statements.

Let $A$ be a $^*$-algebra over $\C$, with or without identity, but with a
non-degenerate product.  The {\it multiplier algebra} $M(A)$ can be characterized as
the largest algebra with identity in which $A$ sits as an essential
two-sided ideal.
We always let $A \otimes A$ denote the algebraic tensor product.
A {\it comultiplication} (or a {\it coproduct}) on $A$ is a non-degenerate $^*$-homomorphism
$\Delta : A \rightarrow M(A \otimes A)$
such that $\Delta (a) (1 \otimes b)$ and $(a \otimes 1)\Delta(b)$ are elements of
$A \otimes A$ for all $a,b \in A$.
It is assumed to be {\it coassociative} in the sense that
$(\Delta \otimes \iota) \circ\Delta(a) = (\iota \otimes \Delta)\circ\Delta (a)$ inside
$ M(A \otimes A\otimes A)$;
where $\iota$ denotes the identity map, see \cite{VD-MH} for a more precise definition.
\begin{defn}
A pair $(A,\Delta)$ of a $^*$-algebra $A$ over $\C$ with a non-degenerate product
and a comultiplication
$\Delta$ on $A$ is called a
{\it multiplier Hopf $^*$-algebra} if the linear maps from
$A \otimes A$  defined by
\begin{align} a \otimes b & \rightarrow \Delta (a)(1 \otimes b) \\
a \otimes b & \rightarrow (a \otimes 1) \Delta (b)
\end{align}
are injective with range equal to $A \otimes A$.
\end{defn}
For any multiplier Hopf $^*$-algebra,
there is a {\it counit} $\varepsilon : A \mapsto \C$
which  is the unique $^*$-homomorphism satisfying
\begin{align}
\label{counit1}& (\varepsilon \otimes \iota)(\Delta (a)(1 \otimes b))
= ab \\
& \label{counit2}(\iota \otimes \varepsilon)((a \otimes 1) \Delta (b)) = ab \end{align}
for all $a,b \in A$.
There is also an {\it antipode} which is the unique linear map
$S : A \mapsto A$ satisfying
\begin{align}
\label{antipode1}
& m(S \otimes \iota)(\Delta (a)(1 \otimes b)) = \varepsilon (a)b \\
 \label{antipode2}& m(\iota \otimes S)((a \otimes 1)\Delta (b)) = \varepsilon (b)a
\end{align}
where $m$ denotes multiplication defined as a linear map from
$A \otimes A$ to $A$.
The antipode is an injective anti-homomorphism and satisfies the relation
$S(a^\ast) = S^{-1}(a)^\ast$
for all $a \in A$.
\begin{defn} A {\it right integral} on $A$ is a linear functional
$ \mu $ \st
\[
(\mu \otimes \iota)(\Delta(a)(1\otimes b))=\mu(a)b.
\]
In general a right integral may not exist, but if it does there is also a left integral (defined in a similar way). This will always be true in the cases we are studying.
A multiplier Hopf $^*$-algebra with a positive right integral is called an {\it algebraic quantum group}, (which should not be confused with a quantization of an algebraic group).
\end{defn}
A multiplier Hopf $^*$-algebra has local units in the following sense:
\begin{lem}
\label{lem:local.units}
If $\A$ is a multiplier Hopf $^*$-algebra and $F$ is a finite subset of $\A$, then
there is
$b\in\A$ \st $ab=a$ for all $a\in F$.
\end{lem}
\begin{proof}
See \cite{DVZ}*{Proposition~2.2}.
\end{proof}

\section{Multiplier Hopf $^*$-algebras in $C_0(G)$}
\label{sec:c-null-G}
We start with $C_0(G)$ where the comultiplication, antipode and counit is  given by
\[
\Delta(f)(x,y)=f(xy)\quad S(f)(x)=f(x\inv)\quad \varepsilon(f)=f(e)
\]
and a Haar integral is given by
$f \mapsto \int fd\mu$
where $\mu$ as before is a left Haar measure on $G$.

\begin{stand}\label{stand:c-null-G}
Let $\A$ be  a  $^*$-subalgebra  of $C_0(G)$ which is also invariant under the antipode $S$.
We also assume that
    \[
    \spn\{\Delta(f)(1\otimes g) \mid f,g\in \A)\}=\A\otimes\A.
\]
 
It can then be shown that $\A$ is
a multiplier Hopf $^*$-algebra with  the coproduct  inherited from $C_0(G)$.
We call $\A$ {\it a multiplier Hopf $^*$-subalgebra}  of $C_0(G)$.
It is actually not necessary to assume that $\A$ is invariant under the antipode $S$,
for details on all this see \cite{DC+VD}.
\end{stand}

The main results in this section is that such a multiplier
Hopf $^*$-subalgebra $\A$ of $C_0(G)$ exists only if $G$ has a compact open
subgroup. If $\A$ also is dense in $C_0(G)$, it is unique and can be described.
\begin{defn}
If $H$ is a compact group, define $\P(H)=$ all polynomial functions on $H$, so
$f\in\P(H)\iff \exists f_i, g_i\in C(H)$ \st
\begin{equation}
\label{defPH}
f(hk)=\sum_1^n f_i(h)g_i(k) \text{ for } h,k \in H.
\end{equation} Thus $\P(H)$ equals the matrix functions corresponding to finite
dimensional (unitary) representations of $H$.
\end{defn}
Suppose $G$ is a locally compact group and that $H$ is a compact open subgroup. Then
$C(H)\subset C_c(G)$ in an obvious way, this is used next.

\begin{lem}
\label{lem:PH} The  following are equivalent\textup:
\begin{enumerate}
\item $f\in \spn\{_x\phi \mid x\in G,\, \phi\in\P(H)\}$.
\item
$\exists f_i\in C_c(G),\, \phi_i\in C(H)$ \st\\
 $f(xh)=\sum_1^n f_i(x)\phi_i(h)$  for $x\in G,\,h \in H$.
\end{enumerate}
\end{lem}
\begin{proof}
(i)$\Rightarrow$(ii): If $\phi\in\P(H)$ satisfies \eqref{defPH}
and  $f= {_x\phi}$, then for $y\in G,\,h \in H$:
\begin{equation}
f(yh)=\Chi_{xH}(yh)\phi(x\inv  yh)= \Chi_{xH}(y)\sum f_i(x\inv
y)g_i(h).
\end{equation}
So $f$ satisfies (ii).

(ii)$\Rightarrow$(i): If $f$ satisfies (ii) and has support in
$\cup_1^m x_iH$ we have $f=\sum \Chi_{x_iH}f$. So we may suppose
$f={_x\phi}$ and have to show that $\phi\in\P(H)$. But this follows
from $\phi(hk)=f(xhk)=\sum f_i(xh)g_i(k)$.
\end{proof}

\begin{lem}
\label{lem:PH2} If $f$ satisfies (ii) above, $f_i, \phi_i$ can be
chosen \st also $f_i$ satisfies (ii) and $\phi_i\in\P(H)$.
\end{lem}
\begin{proof}
First note that we can assume that $\{\phi_i\}$ is an
orthonormal set in $L^2(H)$, so
\begin{align*}
f_i(x)&=\int_H f(xh)\overline{\phi_i(h)}\,dh\\
f_i(xk)&=\int_H f(xkh)\overline{\phi_i(h)}dh =\sum_j f_j(x)\int_H
\phi_j(kh)\overline{\phi_i(h)}\,dh
\end{align*}
therefore $f_i$ satisfies (ii). Since $\spn\{f_i\}$ is $R_H$-invariant
and $\phi_i$ are the matrix functions with respect to this finite
dimensional representation, it follows that $\phi_i\in\P(H)$.
\end{proof}

\begin{lem}
\label{lem:PH=PK} Suppose $H,K$ are two compact open subgroups of
$G$. Then $f$ satisfies the conditions of \lemref{lem:PH} with respect to $H$ if
and only if it does for $K$.
\end{lem}
\begin{proof}
It is enough to show that the statement is true if $H\subset K$.
Using \lemref{lem:PH} it further suffices to show that under the
natural embedding of $C(H)$ into $C(K)$ we have $\P(H)$ mapped
into $\P(K)$. If $K=\cup_1^n k_iH$ we have
\begin{equation}
\Chi_H(kl)=\sum_1^n\Chi_H(kk_i)\Chi_H(k_i\inv l)\righttext{for}
k,l \in K.
\end{equation}
Now suppose $f\in\P(H)$, so $f(hk)=\sum_1^m f_i(h)g_i(k)$  for
$h,k \in H$. Then we have for $k,l \in K$ that
\begin{align*}
f(kl)&=\Chi_H(kl)f(kl)=\sum\Chi_H( kk_i)\Chi_H(k_i\inv l)f(kk_ik_i\inv l)\\
    &=\sum\Chi_H( kk_i)\Chi_H(k_i\inv l)f_j(kk_i)g_j(k_i\inv l)
\end{align*}
which shows that $f\in\P(K)$.
\end{proof}

\begin{defn}
\label{defn:PG} {\it The polynomial functions} on $G$ is
the space $\P(G)$ of all functions $f\in C_c(G)$ satisfying the conditions of \lemref{lem:PH} for some (hence all) compact open subgroups of
$G$.
\end{defn}
\begin{thm}
\label{thm:PG} Suppose $G$ has a compact open subgroup. Then $\P(G)$ is a
multiplier Hopf $^*$-subalgebra
of $C_0(G)$ separating points of $G$ which is invariant under the left and right
action given by $f\mapsto {_xf}$ and $f\mapsto f_x$.
\end{thm}
\begin{proof}
If $f,g\in\P(G)$ there is a compact open subgroup \st the conditions of
\lemref{lem:PH} hold for both. The same subgroup then holds for
both $f+g$ and $fg$.

The antipode in $C_0(G)$ is given by $S(f)(y)=f(y\inv)$. If
$f={_x\phi}$ with $\phi\in\ \P(H)$, then $S(f)={_{x\inv}\psi}$ where
$\psi\in\ \P(xHx\inv)$ is given by $\psi(y)=\phi(x\inv y\inv x)$.
So $\P(G)$ is $S$-invariant. $\P(G)$ is obviously invariant under
$f\mapsto {_xf}$, and if $f$ satisfies \lemref{lem:PH} with respect to  $H$,
then $f_x $ satisfies \lemref{lem:PH} with respect to
$xHx\inv$.

We next have to show that
$\Delta(\P(G))(\P(G)\otimes 1)=\P(G)\otimes\P(G)$ etc.
If $f={_x\phi}$ and $g={_y\psi}$ with
$\phi,\psi\in\P(H)$, then the function
\begin{equation}
 h(s,t):=\Delta(f)(g\otimes 1)(s,t)=\phi(x\inv
st)\psi(y\inv s)
\end{equation}
has support inside $yH\times Hy\inv xH$. By
compactness we get $Hy\inv xH=\cup_1^m h_iy\inv xH$, take $z_i=h_iy\inv
x$ and $K=H\cap_i z_i H z_i\inv$. Suppose
$\Delta(\phi)=\sum_1^n\alpha_j\otimes\beta_j$, then
\begin{align*}
h(s,t)&=\sum_i \Chi_{z_iH}(t)\phi(x\inv st)\psi(y\inv s)\\
      &=\sum_i \Chi_{H}(z_i\inv t)\Chi_H(x\inv st)\phi(x\inv st)\psi(y\inv s)\\
      &=\sum_i \Chi_{H}(z_i\inv t)\Chi_H(x\inv sz_i)\phi(x\inv sz_iz_i\inv
t)\psi(y\inv s)\\
      &=\sum_{i,j} \Chi_{H}(z_i\inv t)\Chi_H(x\inv sz_i)\alpha_j(x\inv sz_i)
        \beta_j(z_i\inv t)\psi(y\inv s).
\end{align*}
From this it follows that $h\in C_0(G)\otimes C_0(G)$.
An easy
computation shows that  $h\in \P(G)\otimes\P(G)$ with respect to the
subgroup $K\times K$.
We have therefore  proved that
$\Delta(\P(G))(\P(G)\otimes 1)\subset\P(G)\otimes\P(G)$. Since the
inverse of the map $a\otimes b\mapsto (a\otimes 1)\Delta(b)$ is
given by
$a\otimes b\mapsto (a\otimes 1)(S\otimes \iota)\Delta(b)$
 it follows that $\Delta(\P(G)(\P(G)\otimes 1)
=\P(G)\otimes\P(G)$. We leave other details to the reader.
\end{proof}
\begin{prop}
\label{prop:PG} $f\in\P(G)$ if and only if $ f\in C_0(G)$ and there are
non-zero functions $g,h, f_i, g_i,f_j', g_j'\in C_0(G)$ \st for
all $x,y\in G$
\begin{enumerate}
\item  $f(xy)g(y)=\sum_1^m f_i(x)g_i(y) $ and
\item
 $f(y)h(xy)=\sum_1^n f'_j(x)g'_j(y) $.
\end{enumerate}
\end{prop}
\begin{proof}
If $f\in\P(G)$ it follows from \eqref{defPH} that (i)
and (ii) hold with $g=h=\Chi_H$. Conversely if
(i) holds, it follows from \corref{cor:mult.C-null-G}  that
there is a compact open subgroup $H$ and functions $h_i\in C_0(G)$
and $k_i\in C(H)$ \st $f(xy)=\sum_1^n h_i(x)k_i(y) $ for $y\in H$.
Finally, from \corref{cor:comp.supp} we have $f\in C_c(G)$, hence
$f\in \P(G)$.
\end{proof}
The next characterization $\P(G)$ will also be useful.
\begin{prop}
\label{prop:PG2} Suppose $G$ has  a compact open subgroup $H$.
Then $f\in\P(G)$ if and only if  there are finitely many functions
$f_i,g_i,f'_j,g_j'\in C_0(G)$ \st
\begin{enumerate}
\item $f(xy)\Chi_H(y)=\sum_1^m f_i(x)g_i(y) $ and
\item
$f(y)\Chi_H(xy)=\sum_1^n f'_j(x)g'_j(y) $.
\end{enumerate}
\end{prop}
\begin{proof} If $f$ satisfies (ii) it follows from \lemref{lem:finite.support} that
$f=\sum_1^n \,_{x_i}\phi_i$ with $\phi_i\in C(H)$ and the sets
$\{x_iH\}$ disjoint. We want to show that (i) implies that
$\phi_i\in \P(H)$:
\begin{equation}
\phi_i(hk)=f(x_i\inv hk)=\sum_1^n f_j(x_i\inv h)g_j(k)
\righttext{for all} h,k \in H,
\end{equation}
so $\phi_i\in \P(H)$ and $f\in\P(G)$ by \lemref{lem:PH}.

Conversely, if $f={_z\phi}$ with $\phi\in\P(H)$ one checks that $f$
satisfies (i) and (ii), so again by \lemref{lem:PH} this is
true for any $f\in\P(G)$.
\end{proof}
\begin{rem}
\label{rem:REM.both} Note that both (i) and (ii) are needed in
general to characterize $\P(G)$: If $G$ is discrete and $H=\{e\}$
then (i) is automatic, if $G$ is compact and $H=G$ then (ii) is
automatic.
\end{rem}
\begin{thm}
\label{thm:A=P(G)} Suppose  $\A$ is a multiplier Hopf
$^*$-subalgebra of $C_0(G)$ separating points. Then $G$ contains a
compact open subgroup $H$ and
 $\A =\P(G)$.
\end{thm}
\begin{proof}
It follows from \corref{cor:mult.C-null-G} and \propref{prop:PG}
that $G$ contains  a compact open subgroup $H$ and that $\A
\subset\P(G)$.

{\it Claim1:} If $\nu$ is a measure on $G$ with compact support
and $f\in\A$, then $f*\nu\in\A$.

Let $C$ be the support of $\nu$. Since $\A$ separates points in
$G$ there is $g\in\A$ \st $g(y)>0$ for $y\in C$. There are
functions $f_i,g_i\in \A$ \st for $y\in C$ we have
\begin{align*}
f(xy\inv)g(y)&=\sum_1^n f_i(x)g_i(y) \\
f(xy\inv)&=\sum_1^n f_i(x)g_i(y)/g(y) \\
f*\nu(x)&=\sum_1^n f_i(x)\nu(g_i/g).
\end{align*}
So $f*\nu\in\A$. In particular this means that $\A$ is invariant
under $f\mapsto f_x$ and therefore also under $f\mapsto {_xf}$.
Moreover it follows that if $f\in\P(G)$ and $g\in\A $ then $f*g\in\A $.

{\it Claim2:} $\Chi_H\in\A$.

By Stone-Weierstrass $\norm{f-\Chi_H}_\infty<\epsilon<1/2$ for some positive function $f\in\A$. Then by Claim1
$g=f*\Chi_H\in\A\cap C_c(G/H)$ and
$\norm{g-\Chi_H}_\infty<\epsilon$. The support of $g$ equals
$\cup^N_{i=0}x_iH$ with $x_0=e$. Take $\alpha_i=g(x_i)$ and define
\begin{equation}
\phi(x)=g(x)\prod_{i=1}^n [\alpha_0 g(x) - \alpha_i g(x\inv x_i)].
\end{equation}
Then $\phi\in\A$, we have $\phi(x_i)=0$ for $i\neq 0$,
$\phi(e)=\alpha_0\prod_1^n [\alpha_0^2  - \alpha_i^2]\neq 0$. So
$\phi=\phi(e)\Chi_H$, hence $\Chi_H\in\A$.

{\it Claim3:} If $f\in\P(H)$ there is $g\in\A$ \st $f=g|_H$.

By taking a minimal decomposition with $f_i, g_i\in C(H)$ \st
    \[
    f(hk)=\sum_1^n f_i(h)g_i(k) \text{ for } h,k \in H
\]
we may assume that $\{g_i\}$ is orthonormal and that $\{f_i\}$ is linearily independent.
Since $\A $ is dense in $C_0(G)$ there are $h_i\in\A $ \st
    \[\int_H g_i(k)h_j(k\inv)\,dk=\delta_{ij}.
\]
Then
$   f*h_i$  is in $\A $ and for $h\in H$
\begin{align*}
f*h_i(h)&=\int_H f(hk)h_i(k\inv)\,dk\\
    &=\int_H\sum_j f_j(h)g_j(k)h_i(k\inv)\,dk
    =f_i(h).
\end{align*}
From this it follows that
    \[
    f(h)=\sum_j f_j(h)g_j(e)=\sum_j f*h_j(h)g_j(e)
\]
which proves the claim.

Finally it follows from Claim2+3 that $ \P(H)\subset\A$,
and then from \lemref{lem:PH}
that $ \P(G)\subset\A$.
\end{proof}

\section{Totally disconnected groups}
\label{sec:tot.disc}

It is natural now to look these groups since they have
a basis of neighborhoods of $e$ consisting of compact open
subgroups. In addition it was our discovery that the {\it smooth}
functions on $G$ is a multiplier Hopf $^*$-algebra that started
this work.

\begin{defn}
\label{defn:smooth}If $G$ is a totally disconnected group, define the {\it smooth}
functions on $G$ by
\begin{align*}
C_c^\infty(G)&=\cup \{C_c(G/H)\mid H \text{ a compact open subgroup} \}\\
             &=\spn \{\Chi_{xH}\mid x\in G,\, H \text{ a compact open subgroup} \}\\
             &=\spn \{\Chi_{xHy}\mid x,y\in G,\, H \text{ a compact open subgroup} \}.
\end{align*}
\end{defn}

\begin{thm}
\label{thm:smooth=PG} If $G$ is a totally disconnected group,
$C_c^\infty(G)=\P(G)$.
\end{thm}
\begin{proof}
If $H$ is a compact open subgroup then $\Chi_H\in\P(G)$, and since both
$C_c^\infty(G)$ and $\P(G)$ are translation invariant
 $C_c^\infty(G)\subset\P(G)$.

To prove the converse, for the same reason it suffices to show
that $\P(H)\subset C_c^\infty(G)$. So suppose $f\in\P(H)$
satisfies
\begin{equation}
 f(hk)=\sum_1^n f_j( h)g_j(k) \righttext{for all} h,k
\in H
\end{equation}
and we may assume that $\{ g_j \}$ is an orthonormal set in
$L^2(H)$. Then as in the proof of \lemref{lem:PH2} we get that
$\{ f_j \}$ is $R_x$-invariant for $x \in H$. This way we get a
finite dimensional representation of $H$ on $X=\spn \{ f_j \}$.
Since $H$ is  totally disconnected,
 by \cite{HR2}*{(28.19)} there is a compact open subgroup $K$ \st
$R_k=I$ on $X$ for $k\in K$. This means that $ f_j $  and
therefore also $ f\in C(H/K) \subset C_c^\infty(G)$.
\end{proof}
\begin{rem}
\label{rem:Bruhat}
Note that if $G$ is totally disconnected  $C_c^\infty(G)$ equals the space of regular functions as defined by  Bruhat in \cite{B}, but in general these spaces are different.
For more about functions on totally disconnected groups, see also
\cite{SilIH}*{chapter 1.1}.
\end{rem}

\section{Multiplier Hopf $^*$-algebras in $C_r^*(G)$}
\label{sec:c*rG}
\begin{defn}
We have already defined the left and right regular representations of $G$ on
$L^2(G)$ in \defnref{defn:reg.repr}. For $f\in L^1(G)$ let
\[
L_f=\int f(x)L_x\, dx\qquad R_f=\int f(x)R_x\, dx.
\]
\end{defn}
Then $C_r^*(G)$ is defined as the norm closure of $\{L_f\mid f\in L^1(G) \}$.
It is standard that $L_x\in M(C_r^*(G))$ and we shall often identify an
element $x\in G$ with $L_x$.
We shall also need the weak closures
\[
\L(G) := \{L_g \mid g \in G\}'' \quad \text{and}\quad \RR(G) := \{R_g \mid g \in
G\}'' .
\]

The comultiplication on $C_r^*(G)$ is defined by
\[
\Delta(L_f)=\int f(x)(L_x\otimes L_x)\, dx
\]
for $f\in L^1(G)$ and  can be extended to a non-degenerate $^*
$-homo\-mor\-phism $C_r^*(G)\mapsto M(C_r^*(G)\otimes
C_r^*(G))$, see \cite{V}*{Proposition~4.3} or (in a more general setting) \cite{L}*{(3.2)}.

The antipode and counit are given by
\[
S(L_f)=\int \Delta_G(x\inv) f(x\inv)L_x\, dx
\qquad
 \varepsilon(L_f)=\int f(x)\, dx,
 \]
where $\Delta_G$ is the modular function of $G$.
A left Haar integral is given by $ w_G(L_f)=f(e)$.

The antipode $S$ can be extended to $C_r^*(G)$, but not the counit $\varepsilon$.
There is an extension of $w_G$  to an (unbounded) weight on $C_r^*(G)$.
For more details we refer to
\cite{GKP}*{Chapter~7.2}.

We shall also need the modular automorphism group corresponding to this weight, it will satisfy
\[
\sigma_t(L_f)=\int \Delta_G(x)^{it} f(x)L_x\, dx .
\]
As usual we also use the notation
\[
\NN_{w_G}=\{a\mid w_G(a^*a)<\infty\}\quad
\MM_{w_G}=\spn\{a^*b\mid a,b\in \NN_{w_G}\}.
\]
\begin{stand}
Let $\A$ be  a  $^*$-subalgebra  of $C_r^*(G)$ which is also invariant under the antipode $S$.
We also here assume that
    \[
    \spn\{\Delta(a)(1\otimes b) \mid a,b\in \A)\}=\A\otimes\A.
\]

It follows that $\A$ is
a multiplier Hopf $^*$-algebra with  the coproduct  inherited from $C_r^*(G)$.
We call $\A$ {\it a multiplier Hopf $^*$-subalgebra}  of $C_r^*(G)$.
As in \ref{stand:c-null-G} it is actually not necessary to assume that $\A$ is invariant under the antipode $S$, for details see \cite{DC+VD}.
\end{stand}

First we address some properties which are not so easy to prove as for $C_0(G)$.
We saw in \secref{sec:prelim} that elements of a multiplier Hopf $^*$-subalgebra
of $C_0(G)$ must have compact support and are therefore automatically integrable
with respect to Haar measure.
We shall see that the similar result is somewhat more complicated in $C_r^*(G)$.
\begin{prop}
\label{prop:analytic} Let $\A$ be a multiplier Hopf $^*$-subalgebra  of
$C_r^*(G)$.
Then $\A$ is $\sigma$-invariant and every element $a\in \A$ is analytic with respect
to the modular automorphism group $\sigma_t$ of the weight $w_G$.
\end{prop}

\begin{proof}
For $a,b\in \A$ we have elements $a_i,b_i\in \A$ \st
\begin{equation}
a\otimes b=\sum_1^n \Delta(a_i)(1\otimes b_i).
\end{equation}
Since $\sigma_t(L_x)=\Delta_G(x)^{it}L_x$,
we have  $(\sigma_t\otimes\sigma_{-t})\circ\Delta=\Delta$ and
\begin{equation}
\sigma_{t}(a)\otimes \sigma_{-t}(b)=\sum_1^n \Delta(a_i)(1\otimes \sigma_{-t}(b_i)).
\end{equation}
Multiply with $1\otimes b^*$ to get
\begin{equation}
\sigma_{t}(a)\otimes b^*\sigma_{-t}(b)=\sum_1^n (1\otimes b^*)\Delta(a_i)(1\otimes
\sigma_{-t}(b_i)).
\end{equation}
Since $a_i,b_i\in \A$, we have $(1\otimes b^*)\Delta(a_i)=\sum c_{ij}\otimes
d_{ij}$, where the sum is finite and the set $\{c_{ij}\}$ is linearily independent.
Take $V_0=\spn\{c_{ij}\}$, then $(1\otimes b^*)\Delta(a_i)\in V_0\otimes\A$ and also
$\sigma_{t}(a)\otimes b^*\sigma_{-t}(b)\in V_0\otimes C_r^*(G)$.
With $b\neq 0$ we see that there is $\epsilon>0$ \st $\sigma_{t}(a)\in V_0$ for
$|t|<\epsilon$. From part (i) of \lemref{lem:mult.Uopen} we see that
$\sigma_{t}(c_{ij})\in V_0$ for all $t\in \R$.

Take
\begin{align*}
V_1&=\spn\{\sigma_{t}(a)\,|\, t\in \R\} \quad\text{and } \\
V_2&=\spn\{\int_{-\infty}^\infty e^{-k(t-t_0)^2}\sigma_{t}(a)\,dt\,|\, t_0\in\R,\,k>0\}.
\end{align*}
If $\alpha$ is a linear functional on $V_1$ which is zero on $V_2$,
we have
\begin{equation}
\int_{-\infty}^\infty e^{-k(t-t_0)^2}\alpha(\sigma_{t}(a))\,dt=0
\end{equation} for all $t_0$ and $k>0$.
This is only possible if
$\alpha(\sigma_{t}(a))\equiv 0$, so  $\alpha= 0$ on $V_1$.
It follows that $V_1=V_2$, and from
\cite{T2}*{Lemma~2.3} that
 $a$ is analytic with respect to the modular automorphism group $\sigma_t$.
 \end{proof}

\begin{prop}
\label{prop:Integrable.C*} Let $\A$ be a multiplier Hopf $^*$-subalgebra  of
$C_r^*(G)$. Then all elements of $\A$
are integrable with respect to the modular automorphism group $\sigma_t$ and
 $\A\subset \NN_{w_G}\cap\MM_{w_G}$.
\end{prop}
\begin{proof}
Take $a\in\A$, we just proved that there is a finite dimensional subspace $V_0$  \st
$\sigma_{t}(a)\in V_0$ for all $t$.
By \lemref{lem:local.units} there is $e\in\A$ \st $ex=x$ for all $x\in V_0$.
Now take $z\in \NN_{w_G}$ \st
$\Norm{e-z}< (4\Norm{e}+2)\inv$ and $y=z^*z$. Then $\Norm{e^*e-y}< \frac{1}{2}$ and
\begin{equation}
a^*a=a^*(e^*e-y)a+a^*ya \leq \frac{1}{2} a^*a+a^*ya,
\end{equation}
so $a^*a\leq 2 a^*ya$. Since $a$ is analytic  with respect to $w_G$,
it follows from \cite{T2}*{Lemma~2.4} that $w_G(a^*ya)<\infty $.
So $w_G(a^*a)<\infty $ and $a \in \NN_{w_G}$.
Since $\A^2=\A$ (as remarked above),
 we also have $a \in \MM_{w_G}$.
 \end{proof}
\begin{rem}
\label{rem:Pedersen-ideal}
Actually $\A$ is contained in the Pedersen ideal of $C_r^*(G)$, but is
in general a proper subset. We shall not need this, but the reader may recognize a main
ingredient of \cite{GKP}*{p~175} in the above proof.
\end{rem}

We now come to the first main result about $C_r^*(G)$:
\begin{thm}
\label{thm:existence.C*r(G)} Suppose $C_r^*(G)$ contains  a multiplier
Hopf $^*$-subalgebra $\A$.
Then $G$ has a compact open subgroup and every element of $\A$ is of the form
$L_\phi$ with $\phi\in C_c(G)$.
\end{thm}
\begin{proof}
By assumption we have $a,b,a_i,b_i\in\A$ with
\begin{equation}
\label{basic.mult}
\Delta(a)(1\otimes b)=\sum_1^n a_i\otimes b_i\neq 0.
\end{equation}
 Then for all $x,y\in G$ we have
\begin{equation}
(1\otimes yx\inv)\Delta(xa)(1\otimes b)=\sum_1^n xa_i\otimes yb_i.
\end{equation}
We have $a,b \in \NN_{w_G}$ and there are $\xi,\eta\in L^2(G)$
such that the following expression is not identically zero:
\begin{equation}
\<xa,w_G\>\< yx\inv b\xi\mid\eta\> =\sum_1^n \<xa_i,w_G\>\< yb_i\xi,\mid\eta\>.
\end{equation}
These  functions are in $C_0(G)$ and satisfy \corref{cor:mult.C-null-G}, so
$G$ has a compact open subgroup and the functions are in fact in $C_c(G)$
by \corref{cor:comp.supp}. With $\widehat a(x)=\phi(x\inv a)$ we then have
$a=\int \widehat a(x)L_x\,dx$.
\end{proof}

As in \secref{sec:c-null-G}  we expect that $C_r^*(G)$ has a unique dense
multiplier Hopf $^*$-subalgebra, and this is true:
\begin{thm}
\label{thm:A=C*P(G)} Suppose  $G$ has a compact open subgroup $H$
and that $\A$ is a dense multiplier Hopf $^*$-subalgebra of
$C_r^*(G)$. Then
\[
 \A =\{L_{\phi}\mid \phi\in\P(G)\}.
 \]
\end{thm}
\begin{proof}
We just saw that $\A \subset\{L_{\phi}\mid \phi\in C_c(G)\}.$
Let
$\widehat \A =\{\widehat a \mid a \in \A\}$.
We want to prove that this is a dense multiplier Hopf $^*$-subalgebra of
$C_0(G)$. It follows from our computations in \thmref{thm:existence.C*r(G)}
that
\[
\spn\{\Delta(\widehat a)(1\otimes \widehat b)\mid a,b\in \A\} = \widehat \A\otimes
\widehat \A.
\]
We have $a\otimes b=\sum \Delta(c_i)(1\otimes d_i)$
so
  $\widehat a(x)\widehat b(x)=\sum \widehat c_i(x)\widehat d_i(e)$,
and therefore $\widehat \A$ is an algebra under pointwise multiplication.
With $b=S(a^*)$ we have $\widehat b(x)=\overline{\widehat a(x)}$, so $\widehat \A$
is conjugation invariant.
By repeating such computations in
various forms, the reader should be convinced that
$\{\widehat a\mid a\in\A\}$ is a multiplier Hopf $^*$-subalgebra of $C_0(G)$.
The conclusion now follows from \thmref{thm:A=P(G)}.
\end{proof}

\begin{rem} In the last part of this section we show that if $G$ has a compact open subgroup $H$,
the unique dense multiplier Hopf $^*$-subalgebra  $\A$ of
$C_r^*(G)$ can be characterized using the  conditional expectation $E:C_r^*(G)\mapsto C_r^*(H)$.
We believe this is useful for generalizations.
\end{rem}

Next we shall give an alternate description of $\A$  which is the dual of \propref{prop:PG2}.
Two tools are needed: the projection
\begin{equation}
\label{pH} p_H=\int_H L_h dh
\end{equation}
(we assume the Haar measure is normalized  such that $\mu(H)=1$)
and the conditional expectation $E:C_r^*(G)\mapsto C_r^*(H)$ given
by
\[
\label{defE}
         E(a)=(\iota\otimes\tau)\Delta(a) =(\tau\otimes \iota)\Delta(a)
\]
where $\tau$  is the vector state given by
$\tau(a)=\<a\Chi_H,\Chi_H\>$. Note that
\[
\Delta\circ E=(E\otimes i)\circ\Delta=(i\otimes E)\circ\Delta
\]
and that for $b\in C_r^*(H)$:
\[
bp_H=\tau(b)p_H,\quad \Delta(b)(1\otimes p_H)=b\otimes p_H.
\]

\begin{lem}
\label{lem:pH1} Suppose $a, a_i, b_i\in C_r^*(G)$  satisfy
\begin{equation}
\label{pH1} \Delta(a)(1\otimes p_H)=\sum_1^n a_i\otimes b_i .
\end{equation}
Then there is a finite set $F$ \st $E(x\inv a)=0$ for $x\notin
FH$ and $a=\sum_{x\in F} xE(x\inv a)$.
\end{lem}
\begin{proof}
By multiplying \eqref{pH1} to the left with $x\inv\otimes y$
and applying $E\otimes\tau$ we get
\[
\Chi_H(yx)E(x\inv a)=\sum \tau( y b_i)E(x\inv a_i) .
\]
Now \lemref{lem:finite.support} gives a finite set $F$ \st
$E(x\inv a)=0$ for $x\notin FH$.

To prove the last claim, choose $F$ \st $FH=\cup_{x\in F} xH$ is a disjoint union
and take $b=\sum_{x\in F} xE(x\inv a)$.
Then $E(y\inv b)=E(y\inv a)$ for all $y\in G$
(look at $y\in FH$ and $y\notin FH$ separately).

So $E(cy\inv b)=E(cy\inv a)$ for all $y\in G$ and $c\in C_r^*(H)$.
Since $\cup y C_r^*(H)$ is dense in $C_r^*(G)$ and $E$ is faithful it follows that $a=b$.
\end{proof}

\begin{lem}
\label{lem:pH2} Suppose $a=\sum_1^n x_i a_i$ with $a_i\in
C_r^*(H)$, $x_i\in G$, $x_j\inv x_i\notin H$ for $i\neq j$ and
that
\begin{equation}
\label{pH2} (a\otimes 1)\Delta(p_H) =\sum_1^m b_k\otimes c_k .
\end{equation}
Then also each $a_i$ satisfies \eqref{pH2}, in fact
\begin{equation}
\label{} (a_i\otimes 1)\Delta(p_H) =\sum_k E(x_i\inv
b_k)\otimes c_k .
\end{equation}
\end{lem}
\begin{proof} Just use the map $b\otimes c\mapsto E(x_i\inv b)\otimes c$ on
\begin{equation}
\label{} \sum (x_i a_i\otimes 1)\Delta(p_H)
=\sum b_k\otimes c_k .
\end{equation}
\end{proof}

\begin{lem}
\label{lem:C*(H)} Suppose $H$ is compact and that $a, b_i,
c_i\in C^*_r(H)$ satisfies
\begin{equation}
\label{C*(H)} (a\otimes 1)\Delta(p_H) =\sum_1^n b_i\otimes c_i.
\end{equation}
Then there is $f\in\P(H)$ \st $a=L_f$.
\end{lem}
\begin{proof}
We may assume that $\{c_i\}$ is linearly independent, so there is
a central projection $e_0\in C^*_r(H)$ \st also $\{c_ie_0\}$ is
linearly independent.
Choose $\psi_j\in (C^*_r(H)e_0)^*$ \st
$\psi_j(c_ie_0)=\delta_{ij}$ and note that $\psi_j$ can be
considered an element of $\P(H)$. Use $i\otimes\psi_j$ on
\eqref{C*(H)} to obtain $b_j=aL_{\psi_j}\in C^*_r(H)S(e_0)$.
So $b_i=L_{f_i}$ for some $f_i\in\P(H)$,
and $(a\otimes 1)\Delta(p_H) =\sum_1^n L_{f_i}\otimes c_i$.
By \eqref{antipode2} $a= \sum_1^n L_{f_i} S( c_i)$ and since
$\{L_{f}\mid f\in\P(H)\}$ is an ideal in $C^*_r(H)$, we have
$a\in\{L_{f}\mid f\in\P(H)\}$.
\end{proof}
\begin{thm}
\label{thm:pC*r(G)} If $G$ has a compact open subgroup $H$ and
$a\in C^*_r(G)$ the following are equivalent:
\begin{enumerate}
    \item
    $a=L_{\phi}$ with $ \phi\in\P(G)$
    \item
    There are finitely many $b_i,c_i,b_j',c_j'\in C^*_r(G)$ \st
    \begin{equation*}
\Delta(a)(1\otimes p_H)=\sum b_i\otimes c_i \quad\text{and}\quad
(a\otimes 1)\Delta(p_H)=\sum b_j'\otimes c_j'.
\end{equation*}
\end{enumerate}
\end{thm}
\begin{proof}
That (i) implies (ii) is left to the reader   (use
\propref{prop:PG2}, multiply with $L_x\otimes L_y$ and integrate).
Conversely, if $a$ satisfies
(ii) it follows from the previous that $a=\sum_k L_{x_k}a_k$ with
$x_k\in G$ and $a_k\in C^*_r(H)$.
 By \lemref{lem:PH2} and \lemref{lem:C*(H)} we get $\phi_k\in\P(H)$ \st
$a_k=L_{\phi_k}$,
so  $f=\sum_k  {_{x_k}{\phi_k}}$
is in $\P(G)$
and we have $a=L_{f}$.
\end{proof}
\begin{rem}
\label{rem:REM.both2} Note that as in \remref{rem:REM.both} both parts of (ii) are
needed in
general to characterize $\A$.
\end{rem}

\section{Multiplier Hopf $^*$-algebras in $C^*(G)$}
\label{sec:c*G}

What happens if we look at $C^*(G)$ instead of $C_r^*(G)$?
Here $C^*(G)$ is the enveloping
$C^*$-algebra of $L^1(G)$ and  the maps $\Delta,\, S,\, \varepsilon,\,
\sigma_t$ in \secref{sec:c*rG} all extends to
$C^*(G)$, \cf\  \cite{iorio}*{Theorem~3.9} or \cite{rae:Roma} for an updated survey.
If $\pi_r$ is the natural map $C^*(G)\to C_r^*(G)$, we also get a weight on
$C^*(G)$ by $a\mapsto \phi(\pi_r(a))$, but this weight is in general
not faithful so $C^*(G)$ is not really a locally compact quantum group.

In \thmref{thm:A=C*P(G)} we showed that the existence of one finite set of elements in
$C_r^*(G)$ satisfying \eqref{basic.mult} implies the existence of a compact open subgroup.
However, this is not true for $C^*(G)$.
Akemann and Walter proved (see \cite{AW} or \cite{Valette}) that if $G$ has property (T),
then there is a  central minimal projection $p_0\in C^*(G)$ \st $\pi_0(p_0)=1$ for the trivial representation $\pi_0$  and $\pi(p_0)=0$ for all other irreducible representations of $G$.
Clearly
\[
\Delta(p_0)(1\otimes p_0)= p_0\otimes p_0,
\]
but there are groups with property (T) -- \eg\ $\SL(3,\R)$ -- which do not have compact open subgroups.

Note that  if $G$ has a compact open subgroup $H$ the analogue of
\thmref{thm:pC*r(G)} can be proved the same way,
 since by \cite{rie:induced}*{Proposition~1.2} there is a conditional expectation
$E:C^*(G)\mapsto C^*(H)=C^*_r(H)$. The map $\tau$ is then defined by
$\tau(a)=\<E(a)\Chi_H,\Chi_H\>$, the proof of \thmref{thm:pC*r(G)} can be repeated verbatim and we have:
\begin{thm}
\label{thm:pC*(G)} If $G$ has a compact open subgroup $H$,
$a\in C^*(G)$ and $U$ is the universal representation of $G$ the following are equivalent:
\begin{enumerate}
    \item
    $a=U_{f}$ with $ f\in\P(G)$
    \item
    There are finitely many $b_i,c_i,b_j',c_j'\in C^*(G)$ \st
    \begin{equation*}
 \Delta(a)(1\otimes p_H)=\sum b_i\otimes c_i \quad\text{and}\quad
(a\otimes 1)\Delta(p_H)=\sum b_j'\otimes c_j'.
\end{equation*}
\end{enumerate}
\end{thm}

\section{Multiplication and convolution operators}
\label{sec:mult+conv}

The dual  locally compact quantum groups $C_0(G)$ and $C_r^*(G)$ have both natural representations on
$L^2(G)$ and we shall study properties of these representations which also turns out to be
equivalent to the existence of a compact open subgroup.
  It is well known, see \cite{Takai}*{Proposition~3.3} or \cite{V}*{Lemme~5.2.8} (although the result is probably older)
that if $a\in C^*_r(G)$ and $f\in C_0(G)$,
then $aM(f)$ is a compact operator on $L^2(G)$).
(See also \cite{BuSm} for a  study of multiplication and convolution operators over $L^p(G)$).

In this section we shall see that $aM(f)$ can not be non-zero and of
\emph{finite rank} unless $G$ has a compact open subgroup.
We shall also see that $aM(f)=M(f)a\neq 0$ is possible
only if $G$ has a compact open subgroup.
We first need the following two results:
\begin{thm}
\label{thm:H.is.comp.open} For a closed subgroup $H$ of $G$,
\begin{enumerate}
\item
$C_0(G)\cap L^\infty(G/H)=C_0(G/H)$ if $H$ is compact and trivial otherwise.
 \item
 $C^*_r(G)\cap \L(H)=C^*_r(G)\cap L^\infty(H\bsl G)' =C^*_r(H)$ if $H$ is open and trivial otherwise.
\end{enumerate}
\end{thm}
\begin{proof} The first statement is  obvious. It follows from the Takesaki-Nielsen-Rieffel
commutant theorem \cite{rie:commut}*{Theorem~2.6} that
\begin{equation} \label{TNR} \L(H)  = \L(G) \cap L^\infty(H \bsl G)'.
\end{equation} Suppose $a \in C_r^*(G)\cap\L(H)$  with $a\gneq 0$. Then there is an
open set $U$ \st $b := M(\Chi_U)aM(\Chi_U)\neq 0$.
So b is a compact operator $L^2(U) \mapsto L^2(U)$
and by the spectral theorem there is
$\lambda\neq 0$ such that the eigenspace
\begin{equation} H_\lambda := \{\xi \mid  b\xi=\lambda\xi\} \end{equation} is finite
dimensional $\neq \{0\}$.
 For $\psi\in L^\infty(H \bsl G)$, $\xi\in H_\lambda$
 then
  \begin{equation}
  bM(\psi)\xi=M(\psi)b\xi=\lambda M(\psi)\xi\,
  \end{equation}  so $M(\psi)H_\lambda \subset H_\lambda$.   We therefore have a
non-zero $\xi\in L^2(U)$
which is an eigenvector for all $M(\psi)$ with $\psi\in L^\infty(H\bsl G)$.
Restricting to $\psi\in C_0(H\bsl G)$
one realizes that there is $x_0 \in G$ (not unique) \st
\begin{equation} \label{tja}
M(\psi)\xi=\psi(x_0)\xi \righttext{for all} \psi\in C_0(H\bsl G).
\end{equation} %
 Let $V = \{x \mid  \xi_0(x)\neq 0\}$, so $\mu(V)>0$  and $\psi(s) = \psi(x_0)$ for
all $s\in V$, $ \psi \in C_0(H\bsl G)$. For this it is necessary that $V\subset
Hx_0$, so  $V V\inv$ is an open subset of $H$ by \cite{HR1}*{(20.17)}
and therefore $H$ is open.
\end{proof}

\begin{thm}
\label{thm:finite.rank}
 Suppose $a \in \RR(G)$
and $f\in L^\infty(G)$ \st $M(f)a \neq 0$ has finite rank. Then $G$ has
a compact open subgroup.
\end{thm}
\begin{proof}
Pick a measurable set $C$ with $0<\mu(C)<\infty$ \st $M(\Chi_C f)a \neq 0$,
therefore
we may assume that $f\in L^2(G)$.
Pick $\xi_i, \eta_i\in L^2(G)$ \st
 \[
 \label{finite.rank}
 M(f)a\xi =\sum_1^n \xi_i\<\xi\mid\eta_i\>\righttext{for all }\xi.
 \]

There is $\xi\in C_c(G)$ \st
$M(f)a\xi\neq 0$, using that $aL_x = L_xa$ we get
\begin{align*}
 M(f)aL_x\xi(y) &=\sum_1^n \xi_i(y)\<L_x\xi\mid\eta_i\>
 \righttext{so }\\
f(y)a\xi(x\inv y) &=\sum_1^n \xi_i(y)\<L_x\xi\mid\eta_i\>.
\end{align*}
 The reader should check that
 $x\mapsto \<L_x\xi\mid\eta_i\>$  is in
$L^2(G)$, so by  \corref{cor:mult.C-null-G}  we can conclude that $G$ has
a compact open subgroup.
\end{proof}
\begin{rem}
We clearly have the same result with $a\in
\L(G)$ instead.
\end{rem}
\begin{thm}
\label{thm:commuting} Suppose $a\in C_r^*(G)$ and $f\in C_0(G)$ are
both non-zero \st $aM(f)=M(f)a$. Then $G$ has a compact open
subgroup.
\end{thm}
\begin{proof}
 Fuglede's Theorem \cite{F} implies that $a^*M(f)=M(f)a^*$,  so
\[
\B = \{g \in L^\infty(G) \mid M(g)a = aM(g)\}
 \]
 is a weakly closed
right invariant $^*$-subalgebra of $L^\infty(G) $, so by \cite{TT}*{Theorem 2}
$\B=L^\infty(H\bsl G)$  for some closed subgroup $H$ of $G$. Since
$f$ is a non-zero element of $C_0(G)\cap L^\infty(H\bsl G)$, we
get from part (i) of \thmref{thm:H.is.comp.open} that $H$ is compact.
Since $a$ is a non-zero element of
$C^*_r(G)\cap\L(H)$,  part (ii) of the same theorem gives that $H$ is open.
\end{proof}

The following description may also be useful.
\begin{defn}A non-zero self-adjoint projection $p$ in a multiplier Hopf
$^*$-algebra is called {\it group-like} (\cf\ \cite{LvD} and \cite{TT}*{Theorem 10}) if
\[
\label{gplike} \Delta(p)(p\otimes 1)=\Delta(p)(1\otimes p)= p\otimes p.
\]
\end{defn}
\begin{prop}
\label{prop:gplike} The following are equivalent:
\begin{enumerate}
    \item
    $G$ has a compact open subgroup
    \item
    $C_0(G)$ has a group-like projection
    \item
    $C^*_r(G)$ has a group-like projection.
\end{enumerate}
\end{prop}
\begin{proof} If $G$ has a compact open subgroup $H$, it is easy to check
that $\Chi_H$ is a group-like projection in $C_0(G)$
 and that $p_H=L_{\Chi_H}$
is a group-like projection in $C^*_r(G)$.

If $p$ is a projection in $C_0(G)$, then $p = \Chi_A$ for a
compact open set $A$. It is easy to see that if $p$  is group-like,
then $A$ is a subgroup of G. Finally, it follows from
\cite{TT}*{Section~5} that if $p\in C^*_r(G)$  is group-like, then $p=L_{\Chi_H}$
for some compact open subgroup $H$ of $G$.
\end{proof}

\begin{rem}
Clearly (i-iii) above implies that $C^*(G)$ has a group-like projection.
However, our remarks in \secref{sec:c*G} show that the reverse implication is false.
\end{rem}

\section{Abelian groups}
\label{sec:abelian}

We close with a quick look at abelian groups.
It is a basic fact of classical Fourier analysis that if we have a non-zero function
$f\in C_c(\R^n)$,
 then its Fourier transform $\widehat f$ is analytic and therefore
 does  not have compact support.  For abelian groups in general we have
the following:
\begin{prop} \label{abelian}
If $G$ is abelian, the following are equivalent:
\begin{enumerate}
   \item     $G$ has a compact open subgroup
    \item     There is a non-zero $f\in C_c(G)$ with
    $\widehat f\in C_c(\widehat G)$.
    \end{enumerate}
    \end{prop}
    \begin{proof}
 If $G$ has a compact open subgroup $H$, then  $f=\Chi_H \in C_c(G)$
 and $\widehat f=\Chi_{H^\perp}\in C_c(\widehat G)$; so (i) implies (ii).

The opposite implication will in fact follow from
 \cite{HR1}*{(24.30)}, but we will give a proof that does not depend on the structure theory of locally compact abelian groups.

Suppose there is a non-zero $f\in C_c(G)$ with
$\widehat f\in C_c(\widehat G)$ and that $U$ is a compact neigborhood of $e$.
Then  there is $g\in C_c(G)$ and $\phi\in C_c(\hat G)$ with $g L_yf =L_yf $ and
 $\phi \widehat{L_yf }=\widehat{L_yf }$
for all $y\in U$.

Hence $L_{\widehat\phi}M(g)$ is a compact operator and
$L_{\widehat\phi}M(g)L_yf =L_yf$ for all $y\in U$.
This implies that $\spn\{L_yf \mid y\in U\}$ is finite dimensional and
it follows from \lemref{lem:mult.Uopen} that $G$ has a compact open subgroup.
\end{proof}


\begin{bibdiv}
\begin{biblist}
\bib{AW}{article}{
  author={Akemann, C. A.},
  author={Walter, M. E.},
  title={Unbounded negative definite functions},
  date={1981},
  journal={Can. J. Math.},
  volume={33},
  pages={862\ndash 871},
}
%
\bib{BS}{article}{
  author={Baaj, S.},
  author={Skandalis, G.},
  title={Unitaires multiplicatifs et dualit\'e pour les produits
              crois\'es de {$C\sp *$}-alg\`ebres},
  date={1993},
  journal={Ann. Sci. \'Ecole Norm. Sup. (4)},
  volume={26},
  pages={425\ndash 488},
}
%
\bib{B}{article}{
  author={Bruhat, F.},
  title={Distributions sur un groupe localement compact et applications
              \`a l'\'etude des repr\'esentations des groupes
              {$\wp$}-adiques},
  date={1961},
  journal={Bull. Soc. Math. France},
  volume={89},
  pages={43\ndash 75},
}

\bib{BuSm}{article}{
  author={Busby, R. C. },
  author={Smith, H. A.},
  title={Product-convolution operators and mixed-norm spaces},
  date={1981},
  journal={Trans. Amer. Math. Soc.},
  volume={263},
  pages={309\ndash 341},
}

\bib{DC+VD}{article}{
   author={De Commer, K.},
  author={Van~Daele, A. },
    title={Multiplier Hopf algebras imbedded in C$^*$-algebraic quantum groups},
  date={2006},
  journal={preprint, /math.OA/0611872},
  volume={},
  pages={1--24},
}

\bib{DVZ}{article}{
  author={Drabant, B. },
  author={Van Daele, A. },
  author={Zhang, Y.},
  title={Actions of multiplier {H}opf algebras},
  date={1999},
  journal={Comm. Algebra},
  volume={27},
  pages={4117\ndash 4172},
}

\bib{F}{article}{
  author={Fuglede, B.},
 title={A commutativity theorem for normal operators},
  date={1950},
  journal={Proc. Nat. Acad. Sci. U. S. A.},
  volume={36},
  pages={35\ndash 40},
  }

\bib{HR1}{book}{
  author={Hewitt, E. },
  author={Ross, K. A.},
  title={Abstract harmonic analysis. {V}ol. {I}: {S}tructure of
              topological groups. {I}ntegration theory, group
              representations},
  publisher={Academic Press Inc.},
  address={Publishers, New York},
  date={1963},
}
\bib{HR2}{book}{
  author={Hewitt, E. },
  author={Ross, K. A.},
  title={Abstract harmonic analysis. {V}ol. {II}: {S}tructure and
              analysis for compact groups. {A}nalysis on locally compact
              {A}belian groups},
  publisher={Springer-Verlag},
  address={New York},
  date={1970},
}

\bib{iorio}{article}{
  author={I{\'o}rio, V. B. M.},
  title={Hopf {$C\sp{\ast} $}-algebras and locally compact groups},
  date={1980},
  journal={Pacific J. Math.},
  volume={87},
  pages={75\ndash 96},
}
\bib{KLQ}{article}{
  author={Kaliszewski, S.},
  author={Landstad, M. B.},
  author={Quigg, J.},
  title={Hecke $C^*$-algebras, Schlich\-ting completions,
and Morita equivalence},
  date={},
  journal={Proc. Edinburgh Math. Soc.},
  pages={(to appear)},
}


\bib{L}{article}{
  author={Landstad, M. B.},
  title={Duality theory for covariant systems},
  date={1979},
  journal={Trans. Amer. Math. Soc.},
  volume={248},
  pages={223\ndash 267},
}
\bib{LvD}{article}{
  author={Landstad, M. B.},
    author={Van~Daele, A.},
title={Compact and discrete subgroups of algebraic quantum groups I},
  date={2007},
  journal={preprint, /math.OA/0702458},
  volume={},
  pages={1\ndash 47},
}

\bib{GKP}{book}{
  author={Pedersen, G.~K.},
  title={{$C\sp{\ast} $}-algebras and their automorphism groups},
  publisher={Academic Press Inc. },
  address={London},
  date={1979},
}

\bib{rae:Roma}{article}{
  author={Raeburn, I.},
  title={Crossed products of {$C\sp *$}-algebras by coactions of
               locally compact groups},
  BOOKTITLE = {Operator algebras and quantum field theory (Rome, 1996)},
      PAGES = {74\ndash 84},
  PUBLISHER = {Internat. Press},
    ADDRESS = {Cambridge, MA},
       YEAR = {1997},
}

\bib{rie:commut}{article}{
  author={Rieffel, M.~A.},
  title={Commutation theorems and generalized commutation relations},
  date={1976},
  journal={Bull. Soc. Math. France},
  volume={104},
  pages={205\ndash 224},
}

\bib{rie:induced}{article}{
  author={Rieffel, M.~A.},
  title={Induced representations of $C^*$-algebras},
  date={1974},
  journal={Adv. Math.},
  volume={13},
  pages={176\ndash 257},
}

\bib{SilIH}{book}{
  author={Silberger, A.~J.},
  title={Introduction to harmonic analysis on reductive $p$-adic groups},
  publisher={Princeton University Press},
  address={Princeton},
  date={1979},
}

\bib{Takai}{article}{
 author={Takai, H.},
  title={On a duality for crossed products of $C^*$-algebras},
date={1975},
  journal={J. Functional Analysis},
  volume={19},
  pages={25\ndash 39},
}
%
\bib{TT}{article}{
  author={Takesaki, M. },
  author={Tatsuuma, N.},
  title={Duality and subgroups},
  date={1971},
  journal={Ann. of Math. (2)},
  volume={93},
  pages={344\ndash 364},
}
\bib{T2}{book}{
  author={Takesaki, M.},
  title={Theory of operator algebras. {II}},
  publisher={Springer-Verlag},
  address={Berlin},
  date={2003},
}
%
\bib{Valette}{article}{
  author={Valette, A.},
  title={Projections in full {$C\sp *$}-algebras of semisimple Lie groups},
  date={1992},
  journal={Mat. Ann. },
  volume={294},
  pages={277\ndash 287},
}
%

\bib{V}{article}{
  author={Vallin, J.~M.},
  title={$C\sp *$-alg\`ebres de Hopf et $C\sp *$-alg\`ebres de Kac},
  date={1985},
  journal={Proc. London Math. Soc.},
  volume={50},
  pages={131\ndash 174},
}

\bib{VD-MH}{article}{
  author={Van~Daele, A.},
  title={Multiplier Hopf algebras},
  date={1994},
  journal={Trans. Amer. Math. Soc.},
  volume={342},
  pages={917\ndash 932},
}

\bib{VD-Adv}{article}{
  author={Van~Daele, A.},
  title={An algebraic framework for group duality},
  date={1998},
  journal={Advances in Mathematics},
  volume={140},
  pages={323\ndash 366},
}

\bib{VDZ}{book}{
  author={Van~Daele, A. },
  author={Zhang, Y.},
    TITLE = {A survey on multiplier {H}opf algebras},
 BOOKTITLE = {Hopf algebras and quantum groups (Brussels, 1998)},
    SERIES = {Lecture Notes in Pure and Appl. Math.},
    VOLUME = {209},
     PAGES = {269\ndash 309},
 PUBLISHER = {Dekker},
   ADDRESS = {New York},
      date = {2000},
}

\end{biblist}
\end{bibdiv}
\end{document}